\documentclass[reqno,11pt]{article}
\usepackage{graphicx, amsmath, amsthm,amssymb}

\input{stochbif.sty}

\def\inst{\mathrm{inst}}
\def\ms{\mathrm{s}}
\def\mt{\mathrm{t}}

\begin{document}


\title{Anomalous behavior of the Kramers rate \\ at bifurcations 
in classical field theories}
\author{Nils Berglund and Barbara Gentz}
\date{Revised version, December 10, 2008}

\maketitle

\begin{abstract}
We consider a Ginzburg--Landau partial differential equation in a bounded
interval, perturbed by weak spatio--temporal noise. As the interval length
increases, a transition between activation regimes occurs, in which the
classical Kramers rate diverges~\cite{Maier_Stein_PRL_01}. We determine a
corrected Kramers formula at the transition point, yielding a finite, though
noise-dependent rate prefactor, confirming a conjecture by Maier and
Stein~\cite{Maier_Stein_SPIE_2003}. For both periodic and Neumann boundary
conditions, we obtain explicit expressions for the prefactor in terms of Bessel
and error functions.
\end{abstract}




\section{Introduction}

Weak noise acting on spatially extended systems can cause a wide range
of interesting phenomena. In particular, it can induce rare transitions between
states which would be otherwise invariant, e.g. nucleation of one phase within
another~\cite{Langer67}, micromagnetic domain
reversal~\cite{Neel49,Braun_PRL93,BrownNovotnyRikvold00}, pattern nucleation in
electroconvection~\cite{Cross_Hohenberg_RMP93}, instabilities in metallic
nanowires~\cite{BurkiStaffordStein_PRL_05}, and many others. The rate of such
transitions for weak noise intensity $\eps$ is in general governed by Kramers'
law\footnote{Throughout this Letter, the notation $a\simeq b$ indicates that
$\lim_{\eps\to0}a/b=1$.} $\Gamma\simeq\Gamma_0\exp\set{-\Delta W/\eps}$, where
the activation energy $\Delta W$ is the energy difference between stable and
transition states, and the rate prefactor $\Gamma_0$ is related to second
derivatives of the system's energy functional at these
states~\cite{Eyring,Kramers}. 

In a series of recent
works~\cite{Maier_Stein_PRL_01,Maier_Stein_SPIE_2003,Stein_JSP_04}, Maier and
Stein studied transition rates in a Ginzburg--Landau partial differential
equation on a finite interval, perturbed by space-time white noise. They
discovered the striking fact that as the interval length approaches a critical
value, which depends on the boundary conditions (b.c.), the rate prefactor
$\Gamma_0$ \emph{diverges}. Although this divergence is reminiscent of the
behavior of certain thermodynamic quantities at phase transitions, it has a
different origin~\cite{Stein04}: It is due to the fact that the Kramers law only
takes into account the effect of quadratic terms in the energy functional on
thermal fluctuations, while at the critical length some quadratic terms vanish
due to a bifurcation, and higher-order terms come into play. 

Maier and Stein conjectured~\cite{Maier_Stein_PRL_01} that the actual rate
prefactor at the bifurcation point behaves like $\Gamma_0\simeq
C\eps^{-\alpha}$, for some constants $C,\alpha>0$. Until recently, 
Kramers rate theory was not sufficiently sharp to allow for the computation
of these constants. Based on a new approach by Bovier {\it et
al}\/~\cite{BEGK}, we developed a method allowing to compute rate prefactors
for potentials with nonquadratic transition states~\cite{BG08a}. The aim of
this Letter is to illustrate the method by determining the constants $C$ and
$\alpha$ in the case of the Ginzburg--Landau equation.


\section{Model}

Consider a one-dimensional classical field $\phi(x,t)$, subjected to the quartic
double-well potential energy function 
\begin{equation}
 \label{mo1}
V(\phi) = \frac14\phi^4 - \frac12\phi^2\;, 
\end{equation} 
to diffusion and to weak space--time white noise. Its evolution is given by
the stochastic partial differential equation (SPDE) 
\begin{equation}
 \label{mo2}
\sdpar\phi t(x,t) = \sdpar{\phi}{xx}(x,t) + 
\phi(x,t) - \phi(x,t)^3 + \sqrt{2\eps}\,\xi(x,t)\;,
\end{equation} 
where $\xi(x,t)$ denotes space--time Gaussian white noise, i.e., formally,
\begin{equation}
 \label{mo2A}
\expec{\xi(x_1,t_1)\xi(x_2,t_2)}= \delta(x_1-x_2)\delta(t_1-t_2)\;.
\end{equation} 
Here we consider the case of a bounded interval $x\in[0,L]$, and either periodic
or Neumann b.c.\ with zero flux, i.e., $\sdpar\phi x(0,t)=\sdpar\phi x(L,t)=0$. 

Note that $\xi(x,t)$ can be rigorously defined by independent white noises
acting on each Fourier mode. For periodic b.c., this leads to setting 
\begin{equation}
 \label{mo5}
\phi(x,t) = \frac{1}{\sqrt{L}}\sum_{k=-\infty}^\infty \phi_k(t)\e^{2\pi\icx
k x/L}
\end{equation} 
and substituting in~\eqref{mo2}. The resulting system of stochastic differential
equations (SDEs) is given by 
\begin{equation}
 \label{mo6}
\dot\phi_k = -\lambda_k\phi_k - \frac1L \sum_{k_1,k_2,k_3\colon k_1+k_2+k_3=k}
\phi_{k_1}\phi_{k_2}\phi_{k_3} + \sqrt{2\eps}\, \dot{W}^{(k)}_t\;,
\end{equation} 
where $\lambda_k=-1+(2\pi k/L)^2$, and the $W^{(k)}_t$ are by definition
independent Wiener processes (see for instance~\cite{Jetschke_86} for a
discussion of the equivalence of different approaches to SPDEs). In the case of
Neumann b.c., setting 
\begin{equation}
 \label{mo8}
\phi(x,t) 
= \frac{1}{\sqrt{L}}\phi_0(t) + \sqrt{\frac{2}{L}}\sum_{k=1}^\infty
\phi_k(t)\cos(\pi k x/L)\;.
\end{equation} 
yields a similar system of SDEs. 

With the SPDE~\eqref{mo2} we associate the energy functional 
\begin{equation}
 \label{mo3}
\cH[\phi] = \int_0^L \biggbrak{\frac12 (\phi'(x))^2 + V(\phi(x))} \6x\;. 
\end{equation} 
For both periodic and Neumann b.c., the uniform configurations
$\phi_\pm\equiv\pm1$ are stable stationary configurations of the system without
noise. Both are minima of the energy functional, of energy
$\cH[\phi_\pm]=-L/4$.
In terms of the Fourier coefficients, for periodic b.c., the potential energy is
given by 
\begin{equation}
 \label{mo7}
\cH[\phi] = \widehat\cH[\set{\phi_k}] = 
\frac12 \sum_{k=-\infty}^\infty \lambda_k\abs{\phi_k}^2 + 
\frac1{4L} \sum_{k_1+k_2+k_3+k_4=0}\phi_{k_1}\phi_{k_2}\phi_{k_3}\phi_{k_4}\;,
\end{equation} 
and a similar relation can be obtained for Neumann b.c.

The value of the activation energy $\Delta W$ for this model is well
known~\cite{Faris_JonaLasinio82,Maier_Stein_PRL_01}. The Kramers rate prefactor
$\Gamma_0$, however, has only been determined for parameters $L$ in certain
ranges, excluding bifurcation values of the model~\cite{Maier_Stein_SPIE_2003}. 


\section{Transition states and activation energy}

The activation energy is the potential energy difference between the initial
stable state $\phi_-$ and the transition state $\phi_\mt$. The latter is defined
as the configuration of highest energy one cannot avoid reaching, when
continuously deforming $\phi_-$ to $\phi_+$ while keeping the energy as low as
possible. The transition state is a stationary state of the energy functional,
that is, it satisfies $\phi_\mt''(x)=-\phi_\mt(x)+\phi_\mt(x)^3$. In addition,
the Hessian operator $\delta^2\cH/\delta\phi^2$ must have a single negative
eigenvalue at $\phi_\mt$. The corresponding eigenfunction specifies the
direction in which the most probable transition path approaches the transition
state.

The shape of $\phi_\mt$ depends on whether the bifurcation
parameter $L$ is smaller or larger than a critical value, the latter depending
on the chosen b.c.~\cite{Maier_Stein_PRL_01}.


\paragraph{Periodic b.c.}

For $L\leqs2\pi$, the transition state is the identically zero function, which
has energy zero. The activation
barrier has thus value $\Delta W=L/4$. 

For $L>2\pi$, there is a continuous one-parameter family of transition states,
of so-called instanton shape, given in terms of Jacobi's elliptic sine by 
\begin{equation}
 \label{ts1}
\phi_{\inst,\varphi}(x) = 
\sqrt{\frac{2m}{m+1}}\sn\biggpar{\frac{x}{\sqrt{m+1}}+\varphi,m}\;.
\end{equation} 
Here $\varphi$ is an arbitrary phase shift, and $m\in[0,1]$ is a parameter
related to $L$ by 
\begin{equation}
 \label{ts2}
4\sqrt{m+1}\,\JK(m) = L\;,
\end{equation} 
where $\JK(m)$ denotes the complete elliptic integral of the first kind.
Note that $m\to0^+$ as $L$ approaches the critical length $2\pi$ from above.
Computing the energy of any instanton transition
state~\eqref{ts1}, one gets~\cite{Maier_Stein_PRL_01} the activation
barrier 
\begin{equation}
 \label{ts3}
\Delta W = \cH[\phi_\inst] - \cH[\phi_-]
= \frac{1}{3\sqrt{1+m}}
\biggbrak{8\JE(m) - \frac{(1-m)(3m+5)}{1+m}\JK(m)}\;, 
\end{equation} 
where $\JE(m)$ denotes the complete elliptic integral of the second kind.


\paragraph{Neumann b.c.}

In this case, the identically zero solution forms the transition state for
all $L\leqs\pi$, so that the activation barrier has again value $\Delta W=L/4$. 

For $L>\pi$, there are two transition states of instanton shape, given by  
\begin{equation}
 \label{ts4}
\phi_{\inst,\pm}(x) = 
\pm\sqrt{\frac{2m}{m+1}}\sn\biggpar{\frac{x}{\sqrt{m+1}}+\JK(m),m}\;,
\end{equation} 
where the parameter $m\in[0,1]$ is now related to $L$ by 
\begin{equation}
 \label{ts5}
2\sqrt{m+1}\,\JK(m) = L\;. 
\end{equation} 
In this case we have $m\to0^+$ as $L$ approaches the critical length $\pi$ from
above. The activation energy is simply half the activation energy~\eqref{ts3} of
the periodic case~\cite{Maier_Stein_PRL_01}.


\section{Rate prefactor}

The rate prefactor $\Gamma_0$ is usually computed by Kramers'
formula~\cite{Eyring,Kramers}
\begin{equation}
 \label{rp1}
\Gamma_0 \simeq \frac1{2\pi}
\sqrt{\biggabs{\frac{\det\Lambda_\ms}{\det\Lambda_\mt}}}
\abs{\lambda_{\mt,0}}\;.
\end{equation} 
Here $\Lambda_\ms=\tdpar{^2\cH}{\phi^2}[\phi_-]$ denotes the linearized
evolution operator at the stable state $\phi_-$,
$\Lambda_\mt=\tdpar{^2\cH}{\phi^2}[\phi_\mt]$ denotes the linearized evolution
operator at the transition state $\phi_\mt$, and $\lambda_{\mt,0}$ denotes the
single negative eigenvalue of $\Lambda_\mt$.

For instance, for Neumann b.c.\ and $L\leqs\pi$, the eigenvalues of
$\Lambda_\mt=-\6^2/\6x^2+1$ are given by $\lambda_k=-1+(\pi k/L)^2$,
$k=0,1,2\dots$, while the
eigenvalues of $\Lambda_\ms=-\6^2/\6x^2-2$ are of the form $\eta_k=2+(\pi
k/L)^2$. It follows~\cite{Maier_Stein_SPIE_2003} that the rate prefactor
is given by 
\begin{equation}
 \label{rp2}
\Gamma_0 \simeq 
\frac{1}{2\pi} \sqrt{\prod_{k=0}^\infty \frac{2+(\pi k/L)^2}{\abs{-1+(\pi
k/L)^2}}}
=\frac{1}{2^{3/4}\pi} \sqrt{\frac{\sinh(\sqrt{2}L)}{\sin L}}\;. 
\end{equation} 
The striking point is that this prefactor diverges, like $(\pi-L)^{-1/2}$, as
$L$ approaches the critical value $\pi$. In fact, this is due to the Kramers
formula~\eqref{rp1} not being valid in cases of vanishing $\det\Lambda_\mt$. 
To confirm Maier and Stein's conjecture that the rate prefactor at $L=\pi$
behaves like $C\eps^{-\alpha}$ and determine the constants $C$ and $\alpha$, we
have to derive a corrected Kramers formula valid in such cases. This can be
done~\cite{BG08a} by extending a technique initially developed by Bovier
{\it et al}\/~\cite{BEGK}, which we outline now.


\paragraph{Potential theory.}

For simplicity, consider first the case of $d$-dimensional Brownian motion
$W^x_t$, starting in a point $x\in\R^d$. Given a set $A\subset\R^d$, the
expected value $w_A(x)=\E[\tau^x_A]$ of the first time $\tau^x_A$ the Brownian
path hits $A$ is known~\cite{Dynkin65} to satisfy the boundary value problem 
\begin{align}
\nonumber
 \Delta w_A(x) &= 1 & x&\in A^c\;,\\
 w_A(x) &= 0          & x&\in A\;.
 \label{pot1}
\end{align} 
The solution can be written as 
\begin{equation}
 \label{pot2}
w_A(x) = \int_{A^c} G_{A^c}(x,y)\,\6y\;, 
\end{equation} 
where $G_{A^c}$ denotes the associated Green's function, 
satisfying $\Delta_x G_{A^c}(x,y)=\delta(x-y)$ and the b.c. (for instance,
$G_{\R^3}(x,y)=1/(4\pi\norm{x-y})$). Similarly, let
$h_{A,B}(x)=\prob{\tau^x_A<\tau^x_B}$ denote the probability that the Brownian
path starting in $x$ hits the set $A$ before hitting the set $B$. It satisfies
the boundary value problem 
\begin{align}
\nonumber
 \Delta h_{A,B}(x) &= 0 & x&\in (A\cup B)^c\;,\\
\nonumber
 h_{A,B}(x) &= 1          & x&\in A\;, \\
 h_{A,B}(x) &= 0          & x&\in B\;.
 \label{pot3}
\end{align} 
This, however, is also the equation satisfied by the electric potential of a
capacitor, with conductors $A$ and $B$ at respective potential $1$ and $0$. If
$\rho_{A,B}(x)$ denotes the surface charge density on the two conductors, we can
write 
\begin{equation}
 \label{pot4}
h_{A,B}(x) = \int_{\partial A} G_{B^c}(x,y)
\rho_{A,B}(y)\,\6y\;.
\end{equation} 
The capacity of the capacitor is simply the total charge accumulated on one
conductor, divided by the potential difference, which equals one: 
\begin{equation}
 \label{pot5}
\capacity_{A}(B) = \int_{\partial A} \rho_{A,B}(y)\,\6y\;.
\end{equation} 
The key observation is the following. Let $C=\cB_\eps(x)$ be a ball of radius
$\eps$ around $x$, and consider the integral $\int_{\partial C} w_A(z)
\rho_{C,A}(z)\6z$. On one hand, using the expression~\eqref{pot2} of $w_A$,
symmetry of the Green's function and then~\eqref{pot4}, one sees that this
integral is equal to $\int_{A^c}h_{C,A}(y)\6y$. On the other hand, as $w_A$ does
not vary much on the small ball $C$~\cite{BEGK}, we can replace $w_A(z)$ by
$w_A(x)$, and the remaining integral is just the capacity. This yields the
relation
\begin{equation}
 \label{pot7}
\Gamma^{-1} = 
\E[\tau^x_A] = w_A(x) \simeq \frac{\displaystyle\int_{A^c}
h_{\cB_\eps(x),A}(z)\,\6z}{\capacity_{\cB_\eps(x)}(A)}\;. 
\end{equation} 
The interest of this relation lies in the fact that capacities can be estimated
by a variational principle. Indeed, the capacity for unit potential difference
is equal to the total energy of the electric field, 
\begin{equation}
 \label{pot8}
\capacity_A(B) = \int_{(A\cup B)^c} \norm{\nabla h_{A,B}(x)}^2\,\6x
= \inf_h \int_{(A\cup B)^c} \norm{\nabla h(x)}^2\,\6x\;,
\end{equation} 
where the infimum is taken over all twice differentiable functions satisfying
the b.c.\ in~\eqref{pot3}.

If, instead of Brownian motion, we consider the solution of a $d$-dimensional
SDE $\dot x = -\nabla \cH(x) + \sqrt{2\eps}\,\dot W_t$, the above steps
can be repeated, provided we replace $\Delta$ by the generator $\eps\Delta -
\nabla\cH\cdot\nabla$ of the equation (the generator is the adjoint of the
operator appearing in the Fokker--Planck equation). The above relations remain
valid, only with the Lebesgue measure replaced by the invariant measure
$\e^{-\cH(x)/\eps}\6x$. Thus we have 
\begin{equation}
 \label{pot9}
\Gamma = \E[\tau^x_A]^{-1}
\simeq \frac{\capacity_{\cB_\eps(x)}(A)}{\displaystyle\int_{A^c}
h_{\cB_\eps(x),A}(z)\e^{-\cH(z)/\eps}\,\6z}\;,
\end{equation} 
where the capacity can be computed via the Dirichlet form 
\begin{equation}
 \label{pot10}
\capacity_A(B) 
= \inf_h  \eps\int_{(A\cup B)^c} \norm{\nabla h(x)}^2\e^{-\cH(x)/\eps}\,\6x\;.
\end{equation} 
The denominator in~\eqref{pot9} can be easily estimated by saddle-point methods,
using the fact that $h_{\cB_\eps(x),A}$ is essentially $1$ in the basin of
attraction of $x$ and $0$ in the basin of $A$. It is equal to leading order to 
$(2\pi\eps)^{d/2}\e^{-\cH(x)/\eps}/\sqrt{\det(\delta^2\cH/\delta x^2)(x)}$.
A good upper bound of the denominator in~\eqref{pot9} is
obtained by inserting a sufficiently good guess for the potential $h$
in~\eqref{pot10}. Assume, e.g., that near a transition state at $0$, the energy
has the expansion 
\begin{equation}
 \label{pot11}
\cH(x) = -\frac12\abs{\lambda_0} x_0^2 + u(x_1) + \frac12 \sum_{j=2}^{d-1}
\lambda_j x_j^2 + \dots\;,
\end{equation} 
where $u(x_1)$ corresponds to the possibly neutral direction in which a
bifurcation occurs. Choosing $h(x)=f(x_0)$ where $\eps f''(x_0) -
\partial_{x_0}\cH(x_0,0,\dots 0)f'(x_0)$ with appropriate b.c.\ 
and substituting in~\eqref{pot10} yields 
\begin{equation}
 \label{pot12}
\capacity_A(B) \leqs \frac{1}{2\pi}
\sqrt{\frac{(2\pi\eps)^{d-1}\abs{\lambda_0}}{\lambda_2\dots\lambda_{d-1}}}
\int_{-\infty}^\infty \e^{-u(x_1)/\eps}\,\6x_1\;.
\end{equation}  
A matching lower bound for the capacity can be obtained by a slightly more
elaborate argument, see~\cite{BG08a} for details.
If $u(x_1)=\frac12\lambda_1 x_1^2$, the integral has value
$\sqrt{2\pi\eps/\lambda_1}$ and we recover the usual Kramers formula.
However,~\eqref{pot12} applies to other cases as well, e.g. a quartic $u(x_1)$. 

We now return to the SPDE~\eqref{mo2}. We apply the above theory first to a
finite-dimensional approximation of the system~\eqref{mo2}, obtained either by
truncation of high wave numbers in its Fourier transform~\eqref{mo6}, or by
replacing the system by a discrete chain~\cite{BFG06a,BFG06b}, and then taking
the limit. A difficulty is that the error terms will depend on the number of
retained modes (see, for instance,~\cite{Liu_CMS_2003} for  estimates on the
convergence rate of the spectral  approximations). Thus the results below are
for now only formal. The error terms in the capacity can, however, be
controlled~\cite{BG08b}. 


\paragraph{Neumann b.c.}

\begin{figure}
\centerline{\includegraphics*[clip=true,width=80mm]{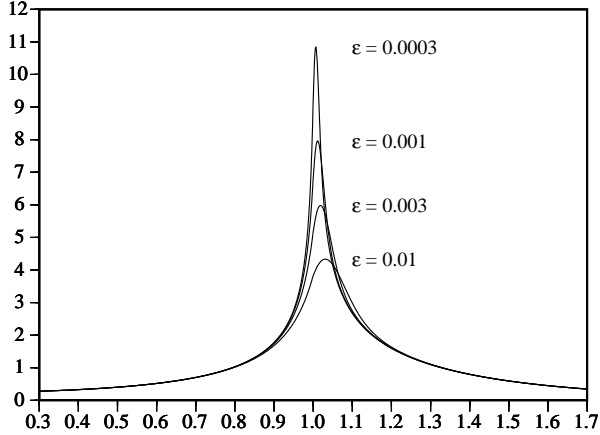}}
\caption{Rate prefactor $\Gamma_0$ as a function of $L/\pi$ for Neumann
 b.c.\ and different values of noise intensity $\eps$.}
\label{fig1}
\end{figure}

The potential energy along the normalized eigenvector in the bifurcating
direction $v_1(x)=\sqrt2 \cos(\pi x/L)$ is
\begin{equation}
 \label{ptn1}
u(\phi_1) = 
\cH[\phi_1v_1] = L \biggbrak{\frac12 \lambda_1\phi_1^2 +
\frac38 \phi_1^4
+ \dots}\;.
\end{equation} 
Evaluating the integral in~\eqref{pot12}, we find~\cite{BG08a} that for
$L\leqs\pi$ the corrected Kramers prefactor to leading order is given by
\begin{equation}
 \label{ptn2}
\Gamma_0 = \frac{1}{2^{3/4}\pi}
\sqrt{\frac{\lambda_1}{\lambda_1+\sqrt{{3\eps/4L}}}}
\Psi_+\biggpar{\frac{\lambda_1}{\sqrt{3\eps/4L}}}
\sqrt{\frac{\sinh(\sqrt{2}L)}{\sin L}}\;, 
\end{equation} 
where $\lambda_1=-1+(\pi/L)^2$ and
$\Psi_+$ is a universal scaling function, given in terms of the modified
Bessel function of the second kind $K_{1/4}$ by 
\begin{equation}
 \label{ptn3}
\Psi_+(\alpha) = \sqrt{\frac{\alpha(1+\alpha)}{8\pi}} \e^{\alpha^2/16}
K_{1/4}\biggpar{\frac{\alpha^2}{16}}\;.
\end{equation} 
For $L\ll\pi$, since $\Psi_+(\alpha)$ tends to $1$ as $\alpha\to\infty$, we
recover to leading order the rate~\eqref{rp2}. For
$\pi-L$ of order $\sqrt{\eps}$, however, the correction terms come into play,
and the factor $\sqrt{\lambda_1}$ in the numerator of~\eqref{ptn2} counteracts
the divergence of the prefactor~\eqref{rp2}. In particular, we have 
\begin{equation}
 \label{ptn5}
\lim_{L\to\pi^-} \Gamma_0 \simeq 
\frac{\Gamma(1/4)}{2(3\pi^7)^{1/4}}
\sqrt{\sinh(\sqrt2\,\pi)}\, \eps^{-1/4}\;. 
\end{equation} 
For $L>\pi$, the rate prefactor is harder to compute, because the transition
states are not uniform. The computation can nevertheless be
done~\cite{Maier_Stein_SPIE_2003} with the help of a method due to
Gel'fand, with the result  
\begin{equation}
 \label{ptn6}
\Gamma_0 \simeq \frac1\pi\bigabs{\mu_0}
\sqrt{\frac{\sinh(\sqrt2\,L)}{\sqrt2\,\abs{(1-m)\JK(m)-(1+m)\JE(m)}}}\;, 
\end{equation} 
where $\mu_0 = 1 - \frac{2}{m+1}\sqrt{m^2-m+1}$ is the negative eigenvalue 
of $\Lambda_\mt$, and $m$ is related to $L$
by~\eqref{ts5}. As $L\to\pi^+$ (that is, $m\to0^+$),
this expression again diverges, namely like $(L-\pi)^{-1/2}$. Proceeding as
above, we find~\cite{BG08a} that the corrected prefactor is obtained by
multiplying~\eqref{ptn6} by 
\begin{equation}
 \label{ptn7}
\frac12 \sqrt{\frac{\mu_1}{\mu_1+\sqrt{3\eps/4L}}} 
\Psi_-\biggpar{\frac{\mu_1}{\sqrt{3\eps/4L}}}\;.
\end{equation} 
Here $\Psi_-$ is again a universal scaling function, given in terms of
modified Bessel functions of the first kind $I_{\pm1/4}$ by 
\begin{equation}
 \label{ptn8}
\Psi_-(\alpha) = \sqrt{\frac{\pi\alpha(1+\alpha)}{32}} \e^{-\alpha^2/64}
\biggbrak{I_{-1/4}\biggpar{\frac{\alpha^2}{64}}
+I_{1/4}\biggpar{\frac{\alpha^2}{64}}}\;, 
\end{equation} 
which converges to $2$ as $\alpha\to\infty$, and $\mu_1$ is the
\emph{second}\/ eigenvalue of $\Lambda_\mt$. We
can in fact avoid the computation of this eigenvalue. Indeed, near the
bifurcation a local analysis shows that
$\mu_1=-2\lambda_1+\Order{\lambda_1^2}=3m+\Order{m^2}$, while further away
from the bifurcation, the quotient in~\eqref{ptn7} is close to $1$. One can
thus replace $\mu_1$ by $3m$ in~\eqref{ptn7}, only causing a multiplicative
error $1+\Order{\eps^{1/4}}$. The resulting behavior of the prefactor
$\Gamma_0$ as $L$ crosses the critical value $\pi$ is shown in Fig.~\ref{fig1}.


\paragraph{Periodic b.c.}

For $L\leqs2\pi$, the transition state is uniform, and the computations are
analogous to those in the previous case. The eigenvalues at the stable and
transition states are now given by $\lambda_k=-1+(2\pi k/L)^2$ and
$\eta_k=2+(2\pi k/L)^2$ with $k\in\Z$, and are thus double except for $k=0$.
This implies that the integral in~\eqref{pot12} is to be replaced
by a double integral over the subspace of the two bifurcating
modes~\cite{BG08b}. The result is 
\begin{equation}
 \label{pp1}
\Gamma_0 \simeq 
\frac1{2\pi} \frac{\lambda_1}{\lambda_1+\sqrt{3\eps/4L}}
\widetilde\Psi_+ \biggpar{\frac{\lambda_1}{\sqrt{3\eps/4L}}} 
\frac{\sinh(L/\sqrt2)}{\sin(L/2)}\;,
\end{equation} 
where the scaling function $\widetilde\Psi_+$ is now given in terms of the
error function by 
\begin{equation}
 \label{pp2}
\widetilde\Psi_+(\alpha) = \sqrt{\frac{\pi}{8}} (1+\alpha) \e^{\alpha^2/8}
\bigbrak{1+\erf(-2^{-3/2}\alpha)}\;.
\end{equation} 
As $\widetilde\Psi_+$ converges to $1$ as $\alpha\to\infty$, for
$2\pi-L\gg\sqrt\eps$, we recover the usual Kramers prefactor, which diverges as
$(2\pi-L)^{-1}$ as $L\to2\pi^-$. However, as $L$ approaches $2\pi$, the
correction terms come into play and we get
\begin{equation}
 \label{pp3}
\lim_{L\to 2\pi^-} \Gamma_0 \simeq
\frac{\sinh(\sqrt2\pi)}{\sqrt3\,\pi}\,\eps^{-1/2}\;. 
\end{equation} 
For $L>2\pi$, we again have to deal with a non-uniform transition
state $\phi_\mt$. An additional difficulty stems from the fact that 
transition states form a continuous family,
so that the Hessian at $\phi_\mt$ always
admits one vanishing eigenvalue. This eigenvalue can be removed by a
regularization procedure due to McKane and Tarlie~\cite{McKane_Tarlie_1995},
which has been applied in
the case of an asymmetric potential in~\cite{Stein_JSP_04}. 
The computations are similar in the symmetric
case~\cite{Stein08_personal}, and yield a rate prefactor per unit length 
\begin{equation}
 \label{pp4}
\frac{\Gamma_0}L \simeq \frac{\abs{\mu_0}}{(2\pi)^{3/2}}  
\sqrt{\frac{2m(1-m)\sinh^2(L/\sqrt2)}{(1+m)^{5/2}\bigabs{\JK(m) -
\frac{1+m}{1-m}\JE(m)}}}  
{\eps}^{-1/2} \;,
\end{equation} 
with $4\sqrt{m+1}\JK(m)=L$ and the same $\mu_0$ as for Neumann b.c. The factor
$\eps^{-1/2}$ reflects the fact that nucleation can occur anywhere
in space~\cite{Stein_JSP_04}. The prefactor now converges to a finite limit
as $L\to2\pi^+$, which differs, however, by a factor $2$ from~\eqref{pp3}. This
apparent discrepancy is solved by applying the corrected Kramers formula, which
shows that~\eqref{pp4} has to be multiplied by a factor 
\begin{equation}
 \label{pp5}
\Phi\biggpar{\frac{3m}{2\sqrt{3\eps/L}}}
\end{equation} 
where $\Phi(x)=\frac12[1+\erf(x/\sqrt2)]$. The resulting rate prefactor is
indeed continuous at $L=2\pi$. 


\section{Conclusion}

We have presented a new method allowing the computation of the Kramers rate
prefactor in situations where the transition state undergoes a bifurcation. In
contrast with the quadratic case, the prefactor is no longer independent of the
noise intensity $\eps$ to leading order, but diverges like $C\eps^{-\alpha}$,
where $\alpha$ is equal to $1/4$ times the number of vanishing eigenvalues. The
constant $C$ can in fact be computed in a full neighborhood of the bifurcation
point, and involves universal functions, depending only on the type of
bifurcation. 
A similar non--Arrhenius behavior of the prefactor has been observed in
irreversible systems~\cite{Maier_Stein_JSP_96}, but there it has an entirely
different origin, namely the development of a caustic singularity in the most
probable exit path.

\paragraph{Acknowledgments.}

We would like to thank Dan Stein for helpful advice, and for
sharing unpublished computations on the periodic-b.c.\ case.
BG was supported by CRC 701 \lq\lq Spectral Structures and Topological Methods
in Mathematics\rq\rq.


\newpage

\small
\bibliography{../BFG}

\def\cprime{$'$}
\providecommand{\bysame}{\leavevmode\hbox to3em{\hrulefill}\thinspace}
\providecommand{\MR}{\relax\ifhmode\unskip\space\fi MR }
\providecommand{\MRhref}[2]{%
  \href{http://www.ams.org/mathscinet-getitem?mr=#1}{#2}
}
\providecommand{\href}[2]{#2}
\begin{thebibliography}{BEGK04}

\bibitem[BEGK04]{BEGK}
Anton Bovier, Michael Eckhoff, V{\'e}ronique Gayrard, and Markus Klein,
  \emph{Metastability in reversible diffusion processes. {I}. {S}harp
  asymptotics for capacities and exit times}, J. Eur. Math. Soc. (JEMS)
  \textbf{6} (2004), no.~4, 399--424.

\bibitem[BFG07a]{BFG06a}
Nils Berglund, Bastien Fernandez, and Barbara Gentz, \emph{Metastability in
  interacting nonlinear stochastic differential equations: {I}. {F}rom weak
  coupling to synchronization}, Nonlinearity \textbf{20} (2007), no.~11,
  2551--2581.

\bibitem[BFG07b]{BFG06b}
\bysame, \emph{Metastability in interacting nonlinear stochastic differential
  equations~{II}: {L}arge-${N}$ behaviour}, Nonlinearity \textbf{20} (2007),
  no.~11, 2583--2614.

\bibitem[BG08a]{BG08b}
Nils Berglund and Barbara Gentz, in preparation, 2008.

\bibitem[BG08b]{BG08a}
\bysame, \emph{The {E}yring--{K}ramers law for potentials with nonquadratic
  saddles}, {\tt arXiv:0807.1681}, 2008.

\bibitem[BNR00]{BrownNovotnyRikvold00}
Gregory Brown, M.~A. Novotny, and Per~Arne Rikvold, \emph{Micromagnetic
  simulations of thermally activated magnetization reversal of nanoscale
  magnets}, vol.~87, AIP, 2000, pp.~4792--4794.

\bibitem[Bra93]{Braun_PRL93}
Hans-Benjamin Braun, \emph{Thermally activated magnetization reversal in
  elongated ferromagnetic particles}, Phys. Rev. Lett. \textbf{71} (1993),
  no.~21, 3557--3560.

\bibitem[BSS05]{BurkiStaffordStein_PRL_05}
J.~Burki, C.~A. Stafford, and D.~L. Stein, \emph{Theory of metastability in
  simple metal nanowires}, Phys. Rev. Lett. \textbf{95} (2005), no.~9, 090601.

\bibitem[CH93]{Cross_Hohenberg_RMP93}
M.~C. Cross and P.~C. Hohenberg, \emph{Pattern formation outside of
  equilibrium}, Rev. Mod. Phys. \textbf{65} (1993), no.~3, 851--1112.

\bibitem[Dyn65]{Dynkin65}
E.~B. Dynkin, \emph{Markov processes. {V}ols. {I}, {II}}, Academic Press Inc.,
  Publishers, New York, 1965.

\bibitem[Eyr35]{Eyring}
H.~Eyring, \emph{The activated complex in chemical reactions}, Journal of
  Chemical Physics \textbf{3} (1935), 107--115.

\bibitem[FJL82]{Faris_JonaLasinio82}
William~G. Faris and Giovanni Jona-Lasinio, \emph{Large fluctuations for a
  nonlinear heat equation with noise}, J. Phys. A \textbf{15} (1982), no.~10,
  3025--3055.

\bibitem[Jet86]{Jetschke_86}
G.~Jetschke, \emph{On the equivalence of different approaches to stochastic
  partial differential equations}, Math. Nachr. \textbf{128} (1986), 315--329.

\bibitem[Kra40]{Kramers}
H.~A. Kramers, \emph{Brownian motion in a field of force and the diffusion
  model of chemical reactions}, Physica \textbf{7} (1940), 284--304.

\bibitem[Lan67]{Langer67}
J.S. Langer, \emph{Theory of the condensation point}, Ann. Phys. \textbf{41}
  (1967), 108--147.

\bibitem[Liu03]{Liu_CMS_2003}
Di~Liu, \emph{Convergence of the spectral method for stochastic
  {G}inzburg-{L}andau equation driven by space-time white noise}, Commun. Math.
  Sci. \textbf{1} (2003), no.~2, 361--375.

\bibitem[MS96]{Maier_Stein_JSP_96}
Robert~S. Maier and D.~L. Stein, \emph{A scaling theory of bifurcations in the
  symmetric weak-noise escape problem}, J. Stat. Phys. \textbf{83} (1996),
  291--357.

\bibitem[MS01]{Maier_Stein_PRL_01}
\bysame, \emph{Droplet nucleation and domain wall motion in a bounded
  interval}, Phys. Rev. Lett. \textbf{87} (2001), 270601--1.

\bibitem[MS03]{Maier_Stein_SPIE_2003}
\bysame, \emph{The effects of weak spatiotemporal noise on a bistable
  one-dimensional system}, Noise in complex systems and stochastic dynamics
  (L.~Schimanski-Geier, D.~Abbott, A.~Neimann, and C.~Van~den Broeck, eds.),
  SPIE Proceedings Series, vol. 5114, 2003, pp.~67--78.

\bibitem[MT95]{McKane_Tarlie_1995}
A.~J. McKane and M.B. Tarlie, \emph{Regularization of functional determinants
  using boundary conditions}, J. Phys. A \textbf{28} (1995), 6931--6942.

\bibitem[N{\'e}e49]{Neel49}
L.~N{\'e}el, \emph{Th\'eorie du trainage magn\'etique des ferro-magn\'etiques
  en grains fins avec application aux terres cuites}, Ann. G\'eophys.
  \textbf{5} (1949), 99--136.

\bibitem[Ste]{Stein08_personal}
D.~L. Stein, private communication.

\bibitem[Ste04]{Stein_JSP_04}
\bysame, \emph{Critical behavior of the {K}ramers escape rate in asymmetric
  classical field theories}, J. Stat. Phys. \textbf{114} (2004), 1537--1556.

\bibitem[Ste05]{Stein04}
\bysame, \emph{Large fluctuations, classical activation, quantum tunneling, and
  phase transitions}, Braz. J. Phys. \textbf{35} (2005), 242--252.

\end{thebibliography}
\bibliographystyle{amsalpha}               

\bigskip\bigskip\noindent
{\small 
Nils Berglund \\ 
Universit\'e d'Orl\'eans, Laboratoire {\sc Mapmo} \\
{\sc CNRS, UMR 6628} \\
F\'ed\'eration Denis Poisson, FR 2964 \\
B\^atiment de Math\'ematiques, B.P. 6759\\
45067~Orl\'eans Cedex 2, France \\
{\it E-mail address: }{\tt nils.berglund@univ-orleans.fr}

\bigskip\noindent
Barbara Gentz \\ 
Faculty of Mathematics, University of Bielefeld \\
P.O. Box 10 01 31, 33501~Bielefeld, Germany \\
{\it E-mail address: }{\tt gentz@math.uni-bielefeld.de}



\end{document}